\documentclass[11pt]{amsart}
\usepackage{epsfig,overpic,comment}
\usepackage[usenames,dvipsnames,svgnames,table]{xcolor}
\usepackage[pagebackref,linktocpage=true,colorlinks=true,linkcolor=Blue,citecolor=BrickRed,urlcolor=RoyalBlue]{hyperref}
\usepackage[msc-links,abbrev,backrefs]{amsrefs} 

\newtheorem{thm}{Theorem}[section]
\newtheorem{lem}[thm]{Lemma}

\newtheorem{prop}[thm]{Proposition}

\theoremstyle{definition}

\newtheorem{note}[thm]{Note}

\newcommand{\R}{\mathbf{R}}

\newcommand{\ol}{\overline}
\newcommand{\C}{\mathcal{C}}

\renewcommand{\d}{\partial}

\renewcommand{\S}{\mathbf{S}}
\renewcommand{\l}{\langle}
\renewcommand{\r}{\rangle}

\renewcommand{\tilde}{\widetilde}

\DeclareMathOperator{\inte}{int}

\title[Rigidity of nonnegatively curved surfaces]{Rigidity of nonnegatively curved surfaces\\relative to a curve}

\author{Mohammad Ghomi}
\address{School of Mathematics, Georgia Institute of Technology,
Atlanta, GA 30332}
\email{ghomi@math.gatech.edu}
\urladdr{www.math.gatech.edu/~ghomi}

\author{Joel Spruck}
\address{Department of Mathematics, Johns Hopkins University,
Baltimore, MD 21218}
\email{js@math.jhu.edu}
\urladdr{www.math.jhu.edu/~js}

\date{\today \,(Last Typeset)}
\thanks{The research of M.G. was supported in part by NSF grant DMS-1711400}

\begin{document}

\begin{abstract}
We prove that any properly oriented $\C^{2,1}$ isometric  immersion of a  positively curved Riemannian surface $M$ into Euclidean 3-space  is uniquely determined, up to a rigid motion, by its values on any curve segment in $M$. A generalization of this result to nonnegatively curved surfaces is presented as well under suitable conditions on their parabolic points. Thus we obtain a local version of Cohn-Vossen's rigidity theorem for convex surfaces subject to a Dirichlet condition. The proof employs in part Hormander's unique continuation principle for elliptic PDEs. Our approach also yields a short proof of Cohn-Vossen's theorem.
\end{abstract}

\maketitle

\section{Introduction}
One of the fundamental results of classical surface theory is Cohn-Vossen's rigidity theorem \cite{cohn-vossen1927, cohn-vossen1936, spivak:v5,sacksteder1962}, which states that isometric closed nonnegatively curved surfaces in Euclidean 3-space are congruent. If the surface is not closed, however,  it generally admits infinitely many noncongruent isometric immersions, and thus other constraints are needed to ensure its rigidity. Here we show that a local Dirichlet condition will suffice. For simplicity, we first state our main result for positively curved surfaces:

\begin{thm}[Main Theorem, First Version]\label{thm:main}
Let $M$ be a connected $2$-manifold and $f$, $\tilde{f}\colon M\to\R^3$ be $\C^{2,1}$ positively curved,   isometric immersions whose mean curvature vectors induce the same orientation on $M$. Suppose that there exists a curve segment $\Gamma$ in $M$ and a proper rigid motion $\rho\colon \R^3\to\R^3$ such that $f=\rho\circ \tilde{f}$ on $\Gamma$. Then $f=\rho\circ \tilde{f}$ on $M$.
\end{thm}

In Section \ref{sec:nonnegative} below we will  generalize the above theorem to nonnegatively curved surfaces, under suitable conditions on their parabolic points. The \emph{manifold} $M$ here may have boundary, and can assume any topological genus \cite{gluck&pan,ghomi&kossowski}. \emph{Isometric} means that  the metrics induced on $M$ by $f$ and $\tilde{f}$ coincide, i.e., 
$
\langle df(v), df(w)\rangle=\langle d\tilde{f}(v), d\tilde{f}(w)\rangle
$ 
for all tangent vectors $v$, $w\in T_p M$ and $p\in M$, where $\langle\cdot,\cdot\rangle$ is the standard inner product in $\R^3$.
If, furthermore, $df(v)\times df(w)$ is parallel to the mean curvature vector of $f$ whenever $d\tilde f(v)\times d\tilde  f(w)$ is parallel to the mean curvature vector of $\tilde f$, we say that the mean curvature vectors  \emph{induce the same orientation} on $M$.
By a \emph{curve segment} in $M$ we mean the image of a smooth embedding $(-\epsilon, \epsilon)\to M$ (which may be arbitrarily small). Finally, a \emph{proper rigid motion}  is an orientation preserving isometry, i.e.,  $\rho\in\mathrm{Iso^+}(\R^3)\simeq\R^3\times \mathrm{SO(3)}$.

The earliest antecedent to Theorem \ref{thm:main}  appears to be a work of John Hewitt Jellett \cite{jellett1849} who in 1849 studied how fixing  a non-asymptotic curve in an analytic surface would render it infinitesimally rigid, see also Weingarten \cite{weingarten}. Later, in 1894, Darboux \cite[Liv. 7, Chap. 5]{darboux:vol3} established  rigidity of analytic surfaces, relative to non-asymptotic curves, via  Cauchy-Kovalevskaya theorem (see Notes \ref{note:CK} and \ref{note:higher}).  Indeed, non-asymptotic curves correspond to non-characteristic hypersurfaces for the underlying PDEs, and fixing a non-asymptotic curve in an isometric embedding  fixes the derivatives of the embedding along that curve (see Note \ref{note:contact}), which furnishes the Cauchy data. These notions are also  implicit in the proofs of Cartan-Janet theorem \cite{spivak:v5,han:book, ivey-landsberg} on analytic isometric embeddings.

As far as we know, Theorem \ref{thm:main} is the first analogue in the smooth category of the Jellett-Darboux rigidity result. Other results relevant to our work include a theorem of Alexandrov and Sen$'$kin \cite{alexandrov-senkin}, also see \cite[p. 181]{pogorelov:book}, who showed that if a pair of isometric positively curved surfaces lie in the upper half-space, are star-shaped and concave with respect to the origin, and  their corresponding boundary points are equidistant from the origin, then they are congruent. There is also a similar result of Pogorelov \cite[p. 178]{pogorelov:book} for convex caps which form concave graphs over the $xy$-plane, and whose corresponding boundary points have equal heights. For more background and references for rigidity problems in surface theory, which date back to Euler, Cauchy, and Maxwell, see \cite{spivak:v5, han:book, pak:book, pogorelov:book,  yau:review, ghomi:problems, sabitov1992}.

The basic outline for proving Theorem \ref{thm:main} is as follows. 
After replacing $\tilde{f}$ with $\rho\circ \tilde{f}$, we assume that $f=\tilde{f}$ on $\Gamma$ and then aim to show that $f=\tilde f$ on $M$. To this end it suffices to establish that $f=\tilde{f}$ on an open neighborhood of a point of $\Gamma$ (Section \ref{sec:local}). This is achieved by showing first  that $f$ and $\tilde{f}$ agree up to second order along $\Gamma$ via geometric arguments (Section \ref{sec:contact}), and then applying a unique continuation principle for elliptic PDEs, with Lipschitz coefficients, due to Hormander (Section \ref{sec:unique}). Finally in Section \ref{sec:nonnegative} we will extend Theorem \ref{thm:main} to the nonnegative curvature case via  works of Sacksteder \cite{sacksteder1962} and Hartman-Nirenberg  \cite{hartman&nirenberg} on parabolic points of surfaces. These methods also yield a short proof of Cohn-Vossen's theorem, which is included in Appendix \ref{appendix:cohn-vossen}.

\begin{note}[Conditions of Theorem \ref{thm:main}]
The orientation condition  in Theorem \ref{thm:main} is necessary. For instance let $M$ be the upper hemisphere of $\S^2$, $f$ be the inclusion map, $\tilde f$ be the reflection of $f$ through the $xy$-plane, and $\Gamma$ be any segment of the boundary  of $M$. Further it is important  that the curvature be positive, at least on $\Gamma$. Consider for instance a flat disk, and roll a corner of it outside its plane. According to \cite[p.212, Rem. 7]{sabitov1992}, there are even negatively curved isometric surfaces which coincide on an open set, but not everywhere else. Thus, in contrast to the Jellett-Darboux result,  Theorem \ref{thm:main} appears to be a strictly elliptic phenomenon. Finally, it is not necessary for $f$, $\tilde f$ to be differentiable everywhere, but it is enough that they be continuous on $M$ while they are $\C^{2,1}$ and isometric on $M\setminus X$, where $X$ is any closed subset without interior points whose complement is connected  and contains $\Gamma$. Then $f=\rho\circ\tilde f$ on $M\setminus X$ and therefore on $M$ by continuity.
\end{note}

\section{Beginning of the Proof: Localization}\label{sec:local}

As mentioned above, we may replace $\tilde{f}$ with $\rho\circ \tilde{f}$ so that $f=\tilde{f}$ on $\Gamma$.   Furthermore,  we may assume that $\Gamma$ lies on the boundary $\d M$ of $M$. Indeed, if a point of $\Gamma$ lies in the interior  of the manifold, $\inte(M):=M\setminus \d M$, we may  extend a small segment containing that point  to a closed curve $\ol \Gamma$ bounding a disk $D\subset\inte(M)$. Then $M':=M\setminus\inte(D)$ and $D$ form a pair of manifolds whose boundaries contain $\Gamma$. So $M$ will be rigid relative to $\Gamma$ if and only if $M'$ and $D$ are rigid relative to $\Gamma$.
To prove Theorem \ref{thm:main}, it suffices then to show that: 

\begin{prop}\label{prop:local}
Let $M$, $f$, $\tilde f$, and $\Gamma$ be as in Theorem \ref{thm:main}.
Suppose that $\Gamma\subset \d M$ and $f=\tilde f$ on $\Gamma$. Then every point of $\Gamma$ has an open neighborhood in $M$ where $f=\tilde{f}$.
\end{prop}

 Indeed, suppose that the above proposition holds, let  $U$ be the union of all open sets in $M$  where $f=\tilde{f}$, and $\d U$ be the topological boundary of $U$ in $M$.  Suppose, towards a contradiction, that there exists a point $q_1\in \d U\cap \inte(M)$, for otherwise we are done. Identify a neighborhood of $q_1$ in $M$ with an open disk $\Omega\subset\R^2$ centered at $q_1$. Let $q_0\in U\cap \Omega$. There exists $\delta>0$ such that the closed disk $B_0$ of radius $\delta$ centered at $q_0$ lies in $U\cap\Omega$. Let $q_t:=(1-t)q_0+t q_1$, and $s\in[0,1]$ be the supremum of $t\in [0,1]$ such that $B_t\subset U$, where $B_t$ is the closed disk of radius $\delta$ centered at $q_t$. Then there exists a point $r\in \d B_s\cap \d U$. Applying Proposition \ref{prop:local} to a segment of $\d B_s$ containing $r$ yields that $f=\tilde{f}$ on an open neighborhood of $r$. Thus $r\not\in \d U$, which is the desired contradiction. So it remains to prove Proposition \ref{prop:local}, which is undertaken in the next two sections.

\section{Order of Contact Along $\Gamma$}\label{sec:contact}
 We say that $f$ and $\tilde f$ have \emph{contact of order} $2$ along $\Gamma$ if, in some local coordinates, their derivatives agree up to second order on $\Gamma$. Here we show that, under the hypothesis of Proposition \ref{prop:local}:

\begin{lem}\label{lem:contact}
 $f$ and $\tilde f$ have contact of order $2$ along $\Gamma$.
\end{lem}

\noindent The above lemma appears to have been known, as a version of  it is discussed in a  Russian text by Kagan \cite[p. 199--200]{kagan:book}.   We include our own treatment here, which will also yield a quick proof of the Jellett-Darboux theorem (Note \ref{note:higher}). 

 To set the stage, we  identify a small neighborhood $\Omega^+$ of a point of $\Gamma$ in $M$ with a half disc in $\R^2$ bordering the $y$-axis, and lying to the right of it. We set $f(t):=f(0,t)$ for any mapping $f$ defined on $\Omega^+$. Then, by assumption,
\begin{equation}\label{eq:f}
f(t)=\tilde f(t).
\end{equation}
 We need to show that the partial derivatives of $f$ and $\tilde{f}$ agree up to second order on $\Gamma$ (the $y$-axis), i.e., $f_i(t)=\tilde{f}_i(t)$ and $f_{ij}(t)=\tilde{f}_{ij}(t)$ for $i$, $j=1$,  $2$.  To this end we first note that
$$
f_2(t)=\tilde{f}_2(t),\quad\quad\text{and}\quad\quad f_{22}(t)=\tilde{f}_{22}(t).
$$
By the isometry assumption we may also record that the coefficients of the induced metric tensor $g$ of $M$ are given by
\begin{equation}\label{eq:gij}
g_{ij}:=\langle f_i, f_j\rangle = \langle \tilde{f}_i, \tilde{f}_j\rangle=:\tilde g_{ij}.
\end{equation}
Further  we may assume that $\{f_1(t), f_2(t)\}$ and $\{\tilde f_1(t),\tilde f_2(t)\}$ are each orthonormal. This may be achieved by letting $\gamma(t)$ denote an arc length parametrization for $\Gamma$ (with respect to $g$), $\nu(t)$ be the inward unit normal vector field along $\Gamma$ (again with respect to $g$),  and resetting 
$$
f(s,t):=f\big(\exp_{\gamma(t)}\big(s\,\nu(t)\big)\big),\quad\quad\text{and}\quad\quad\tilde f(s,t):=\tilde f\big(\exp_{\gamma(t)}\big(s\,\nu(t)\big)\big),
$$
 where $\exp$ is the exponential map of $M$, and $(s,t)$ ranges in a half-disk which we again denote by $\Omega^+$.
Then $f(t)$ has unit speed, and $f_1(t)$, $\tilde f_1(t)$ are inward conormals of $f(t)$ with respect to $f(\Omega^+)$, $\tilde f(\Omega^+)$. More generally,
\begin{equation}\label{eq:fermi}
g_{11}= 1,\quad\quad\text{and}\quad\quad g_{12}= 0
\end{equation}
on $\Omega^+$. These equations hold, via Gauss's Lemma, because $s\mapsto f(s,t)$ traces a geodesic with unit speed. Next note that since $f$ has positive curvature, it has no asymptotic directions. So if 
$$
n(t):=f_1(t)\times f_2(t)\quad\text{and}\quad \tilde n(t):=\tilde f_1(t)\times \tilde f_2(t)
$$
denote the unit normals of $f$ and $\tilde f$ on $\Gamma$, then  $\l f_{22}(t),n(t)\r$ and $\l \tilde f_{22}(t),\tilde n(t)\r=\l f_{22}(t),\tilde n(t)\r$ do not vanish. Further, by the orientation assumption in Theorem \ref{thm:main}, they must have the same sign:
\begin{equation}\label{eq:f22}
\l f_{22}(t),n(t)\r \l f_{22}(t),\tilde n(t)\r>0.
\end{equation}
Indeed, since the curvature is positive, the mean curvature vector points to the side of the tangent plane where the surface locally lies. Thus $\l f_{22}(t),n(t)\r$,  $\l \tilde f_{22}(t),n(t)\r $ are both positive  (negative) if and only if $n$, $\tilde n$ are parallel (antiparallel)  to the mean curvature vectors of $f$, $\tilde f$ respectively.

 \subsection{First order contact} As we already know that $f_2(t)=\tilde f_2(t)$, it remains to check that $f_1(t)=\tilde f_1(t)$. Since $f(t)$ has unit speed, and $f_1(t)$, $\tilde f_1(t)$ are inward conormals, the geodesic curvature of $\Gamma$ with respect to the interior of $M$  is given by
 $$
 \l f_{22}(t), f_1(t)\r= -\frac12 (g_{22})_1(t)=-\frac12 (\tilde g_{22})_1(t)=\l \tilde f_{22}(t), \tilde f_1(t)\r.
 $$ 
By \eqref{eq:f22},  $f_{22}(t)\neq 0$. So the principal normal $N(t):=f_{22}(t)/|f_{22}(t)|$ of $f(t)$   is well defined, and the last displayed expression yields that
 \begin{equation}\label{eq:f1N}
 \l f_1, N\r=\l\tilde f_1, N \r.
\end{equation}
Now if $B(t):=f_2(t)\times N(t)$ denotes the binormal vector of $f(t)$, then $\{N(t), B(t)\}$ forms an orthonormal basis for the normal planes of $f(t)$, which contain $f_1(t)$. Thus
 $$
 \l f_1, N\r^2+\l f_1, B\r^2=|f_1|^2=|\tilde f_1|^2= \l \tilde f_1, N\r^2+\l \tilde f_1, B\r^2.
 $$
 So it follows that $\l f_1, B\r=\pm \l\tilde f_1, B\r$. If  $\l f_1, B\r=-\l\tilde f_1, B\r$, then
 $$
\l n, N\r=\l f_2\times f_1, f_2\times B\r=\l f_1, B\r=-\l\tilde f_1, B\r=-\l f_2\times \tilde f_1, f_2\times B\r=-\l\tilde n, N\r,
$$
which contradicts \eqref{eq:f22}. So we conclude that $\l f_1, B\r=\l\tilde f_1, B\r$ which together with \eqref{eq:f1N} yields that
\begin{equation}\label{eq:f1}
f_1(t)=\tilde f_1(t).
\end{equation}

\subsection{Second order contact}
To show that the second derivatives of $f$ and $\tilde{f}$ match up along $\Gamma$ first note that, since $g_{ij}=\tilde g_{ij}$, 
$$
\langle f_{ij}, f_k\rangle=\Gamma_{ij}^k=\frac12\sum_\ell g^{\ell k}\big((g_{\ell i})_j+ (g_{j\ell})_i-(g_{ij})_\ell\big)=\tilde \Gamma_{ij}^k= \langle \tilde{f}_{ij}, f_k\rangle,
$$
where $\Gamma_{ij}^k$, $\tilde \Gamma_{ij}^k$ are the Christoffel symbols associated to $f$, $\tilde f$, and $(g^{ij}):=(g_{ij})^{-1}$.
So it remains to check that  the coefficients of the second fundamental form $\ell_{ij}:=\langle f_{ij},n\rangle$, $\tilde\ell_{ij}:=\langle \tilde{f}_{ij},n\rangle$ agree on $\Gamma$. To this end note that
$$
\ell_{12}(t)=-\langle f_1(t), n_2(t)\rangle=-\langle \tilde f_1(t),  n_2(t)\rangle=\tilde\ell_{12}(t).
$$
 Further $\ell_{22}(t)=\tilde \ell_{22}(t)$, since $f_{22}(t)=\tilde f_{22}(t)$. By Theorema Egregium, the curvature $K:=\det(\ell_{ij})/\det(g_{ij})$ of $f$ coincides with the curvature $\tilde K:=\det(\tilde\ell_{ij})/\det(\tilde g_{ij})$ of $\tilde{f}$. Indeed, Theorema Egregium does hold for $\C^2$ surfaces \cite{hartman-wintner}.
Thus 
$$
\det(\ell_{ij})=K\det(g_{ij})=\tilde K\det(\tilde g_{ij})=\det(\tilde\ell_{ij}).
$$ 
Furthermore, by \eqref{eq:f22}, $\ell_{22}$, $\tilde \ell_{22}$ do not vanish along $\Gamma$. So
$$
\ell_{11}(t)=\frac{\ell_{12}^2(t)}{\ell_{22}(t)}=\frac{\tilde\ell_{12}^2(t)}{\tilde\ell_{22}(t)}=\tilde\ell_{11}(t).
$$
Hence $\ell_{ij}(t)=\tilde\ell_{ij}(t)$,  as desired, which completes the proof of Lemma \ref{lem:contact}.

\begin{note}[Non-asymptotic curves]\label{note:contact}
In the proof of Lemma \ref{lem:contact}  above we used the positive curvature assumption only to ensure that \eqref{eq:f22} holds. Thus Lemma \ref{lem:contact} holds for any pairs of  $\C^2$ surfaces, regardless of their curvature, as long as $\Gamma$ is \emph{non-asymptotic}, i.e., never tangent to an asymptotic direction, and  the normal curvatures of $f$, $\tilde f$ assume the same sign along $\Gamma$.  
\end{note}

\begin{note}[The analytic case]\label{note:CK}
By differentiating \eqref{eq:gij}, and using \eqref{eq:fermi},
one quickly obtains the following equations, assuming that $f$ is $\C^3$, see  Spivak \cite[p. 150]{spivak:v5} or Han-Hong \cite[p. 6]{han:book}: 
\begin{equation}\label{eq:f11f22}
\langle f_{11},f_{2}\rangle=0,\quad \l f_{11},f_1\r=0,\quad \l f_{11},f_{22}\r=-\frac12(g_{22})_{11}+|f_{12}|^2.
\end{equation}
By Cauchy-Kovalevskaya theorem, these equations have a unique solution once $f(t)$ and $f_1(t)$ have been prescribed, and $f$ is  analytic, see \cite[p. 150--153]{spivak:v5} or \cite[Lem. 1.1.4]{han:book}. Thus, by \eqref{eq:f} and \eqref{eq:f1}, $f=\tilde f$ when $\tilde f$ is also analytic. This proves Theorem \ref{thm:main} in the analytic case, and more generally establishes the rigidity of all analytic surfaces relative to  non-asymptotic curves, as first observed by Darboux \cite[p. 280]{darboux:vol3}; see also Hopf and Samelson \cite[Sec. 3]{hopf-samelson}.   For some isometric extension results in the smooth category see \cite{jacobowitz,ghomi&greene,hungerbuhler-wasem}.
\end{note}

\begin{note}[Higher order contact]\label{note:higher}
When $f$ and $\tilde f$ are $\C^{k}$, it can be shown directly that they coincide up to order $k$ along $\Gamma$. This yields a quick proof of Jellett-Darboux rigidity theorem for analytic surfaces, without invoking the Cauchy-Kovalevskaya theorem. Indeed, we claim that if $f$ and $\tilde f$ agree up to order $2\leq m<k$ on $\Gamma$, then  they agree up to order $m+1$. To see this let $\alpha:=\alpha_1\alpha_2\dots \alpha_m$, where $\alpha_i:=1$, $2$. Then, 
$
f_\alpha(t)=\tilde f_\alpha(t),
$
for all $\alpha$,  
 which yields that $f_{\alpha2}(t)=\tilde f_{\alpha2}(t)$. 
 By commutativity, 
it remains then to check that $f_{\alpha1}(t)=\tilde f_{\alpha1}(t)$, where all $\alpha_i=1$.  Since $\{f_1(t), f_2(t), f_{22}(t)\}$ is linearly independent, due to the non-asymptotic assumption on $\Gamma$, this follows from repeatedly differentiating the equations \eqref{eq:f11f22} with respect to the first variable, which shows that $f_{\alpha1}(t)$ is determined by $f_{\alpha'2}(t)$, $f_{\alpha}(t)$, and  lower order derivatives. Thus, by induction, $f$ and $\tilde f$ agree up to order $k$ on $\Gamma$. Consequently, $f=\tilde f$ on $M$ when $f$ and $\tilde f$ are analytic.
\end{note}

\section{Unique Continuation}\label{sec:unique}
 To complete	 the proof of Proposition \ref{prop:local}, and therefore of Theorem \ref{thm:main}, it remains to show that $f=\tilde f$ on the region $\Omega^+$ discussed in the last section. To this end, let $\Omega^-$ be the reflection of $\Omega^+$ with respect to the $y$-axis, and set $\Omega:=\Omega^+\cup\Omega^-$. We may extend $f$, $\tilde f$ isometrically to all of $\Omega$, without loosing regularity, as follows. By the Lipschitz version of Whitney's extension theorem \cite[Thm. 2.64]{brudnyi:book}, first we extend $f$ to $\Omega$ so that $f\in \C^{2,1}(\Omega,\R^3)$; see also \cite[p. 10]{friedman:book} or \cite[Sec. 5.4]{evans:book} for explicit constructions via ``higher order reflection". Then we extend $\tilde f$ by setting it equal to $f$ on $\Omega^-$.  By Lemma \ref{lem:contact}, $f$ agrees with $\tilde f$ up to order $2$ on the $y$-axis, so  $\tilde f_{ij}$ are continuous on $\Omega$; furthermore, $\tilde f_{ij}$ are Lipschitz both on $\Omega^+$ and $\Omega^-$, which quickly yields that $\tilde f_{ij}$ are Lipschitz on $\Omega$. So $\tilde f\in\C^{2,1}(\Omega,\R^3)$ as well.

After replacing $\Omega$ by a smaller disc, we may assume that $f$ and $\tilde{f}$ are positively curved on $\Omega$.
Now for a unit vector $e\in\R^3$, let $u:=\langle f, e\rangle$, $\tilde u:=\langle \tilde f, e\rangle$. Then $u, \tilde{u}\in \C^{2,1}(\Omega)$, and they both satisfy the Darboux equation \cite[p. 45]{han:book}:
\begin{equation}\label{eq:darboux}
\det (\nabla_{ij}u)=K \det(g_{ij}) (1-|\nabla u|_g^2),
\end{equation}
where $(\nabla_{ij})$ and $\nabla u$ are the Riemannian Hessian and gradient respectively and $|\cdot|_g:=\sqrt{\l\cdot,\cdot\r_g}$ is the Riemannian norm. More explicitly,
$$
\nabla_{ij}u:= u_{ij}-\sum_k\Gamma_{ij}^k u_k,\quad\text{and}\quad \nabla u:=\sum_{ij}g^{ij}u_i\partial_j,
$$
where $\d_j$ denote the standard basis of $\R^2$. Thus \eqref{eq:darboux} is a fully nonlinear Monge-Amp\`{e}re equation of general form  \cite[Sec. 3.8]{trudinger-wang2008}. Since $K>0$, this equation is  elliptic whenever $|\nabla u|_g<1$ \cite[p. 46]{han:book}, which is the case here for an open set of directions $e$. We claim that, for these direction, $\phi:=u-\tilde u$ vanishes identically on $\Omega$, which is all we need. Indeed then we have $u\equiv \tilde u$ for $3$ linearly independent directions, which yields that $f\equiv \tilde f$. To this end, we subtract  \eqref{eq:darboux}  from  its counterpart in terms of $\tilde u$. A straight forward computation yields that
\begin{equation}\label{eq:phi}
\sum_{ij}(\nabla_{ij}^*u+\nabla_{ij}^*\tilde u)\nabla_{ij}\phi+2 K\det(g_{ij})\langle \nabla (u+\tilde u),\nabla\phi\rangle_g=0,
\end{equation}
where $(\nabla_{ij}^*):=\det(\nabla_{ij})(\nabla_{ij})^{-1}$ is the cofactor matrix of $(\nabla_{ij})$, i.e., $\nabla_{11}^*=\nabla_{22}$, $\nabla_{12}^*=-\nabla_{12}$, and 
$\nabla_{22}^*=\nabla_{11}$.
Note that \eqref{eq:phi} is linear in terms of $\phi$, and $\phi$ vanishes on an open subset of $\Omega$.
Hence, by the following unique continuation principle,    see also Armstrong and Silvestre \cite[Prop. 2.1]{armstrong-silvestre} and Garofolo and Lin \cite{garofalo-lin}, $\phi$ vanishes  identically on $\Omega$ as claimed. 

\begin{lem}[Hormander \cite{hormander:vol3}, Thm. 17.2.6]\label{lem:hormander}
Let $\Omega\subset\R^2$ be a connected domain and $u\colon\Omega\to \R$ be a solution to the linear equation
\begin{equation}\label{eq:abc}
\sum_{ij}a_{ij}(x)u_{ij}+\sum_i b_i(x) u_i +c(x) u=0,
\end{equation}
where $a_{ij}\colon\Omega\to\R$ are uniformly elliptic and Lipschitz, while $b_i,c\colon\Omega\to\R$ are bounded measurable functions. If $u$ vanishes on an open subset of $\Omega$, then $u\equiv 0$.
\end{lem}

\emph{Uniformly elliptic} means that  the eigenvalues of $(a_{ij})$ are bounded below by a positive constant. To check this and other requirements needed to apply Lemma \ref{lem:hormander} to \eqref{eq:phi} note that in this context
$$
a_{ij}=\nabla_{ij}^*u+\nabla_{ij}^*\tilde u.
$$
 So $a_{ij}$ are Lipschitz, since $u$, $\tilde u\in\C^{2,1}(\Omega)$. Next note that  $\det(\nabla_{ij}u)$, $\det(\nabla_{ij}\tilde u)>0$ by \eqref{eq:darboux}, since $K>0$ by assumption. Thus eigenvalues of $(\nabla_{ij}u)$ and $(\nabla_{ij}\tilde u)$ never vanish on $\Omega$. Now since, by construction, $u$, $\tilde u$ agree up to second order at some point of $\Omega$,  these eigenvalues will coincide at one point, and thus will always carry the same sign. In particular we may assume that they are positive on $\Omega$, after replacing $e$ with $-e$ if necessary. So $(\nabla_{ij}u)$ and $(\nabla_{ij}\tilde u)$ are positive definite matrices. Consequently their cofactor matrices $(\nabla_{ij}^*u)$ and $(\nabla_{ij}^*\tilde u)$ are positive definite as well. Hence so  is their sum $(a_{ij})$. Now  we may assume that $a_{ij}$ are uniformly elliptic on $\Omega$, after replacing $\Omega$ by a smaller disk with compact closure $\ol\Omega'\subset\Omega$.  Finally note that  $b_i$  are continuous on $\ol\Omega'$, while $c\equiv 0$, so they are all bounded and measurable. This concludes the proof of Theorem \ref{thm:main}.
 
 \section{Nonnegative Curvature}\label{sec:nonnegative}
 
 Here we generalize Theorem \ref{thm:main} to nonnegatively curved surfaces. Note that if the set of parabolic, or zero curvature, points $M^0\subset M$ of a nonnegatively curved $\C^{2,1}$ immersion $f\colon M\to \R^3$ does not have interior points and does not disconnect $M$, then $f=\tilde f$ on $M\setminus M^0$ by Theorem \ref{thm:main}, and therefore, by continuity, $f=\tilde f$ on $M$. Thus the nontrivial  case is when $\inte(M^0)$ is nonempty. 
 
 \begin{thm}[Main Theorem, Full Version]\label{thm:main2}
Let $M$, $f$ and $\tilde f$  be as in Theorem \ref{thm:main}, except that the curvature of $f$ is allowed to be nonnegative. Let $M^0\subset M$ be the set of points where the curvature of $f$ vanishes. Suppose that 
\begin{enumerate}
\item[(i)]{$M^0$ contains no curve $\ell$ such that $f(\ell)$ is a complete line,}
\item[(ii)]{$M^0$ is complete, i.e., its Cauchy sequences converge,}
\item[(iii)]{$M^0\subset\inte(M)$,}
\item[(iv)]{$M\setminus M^0$ is connected.}
\end{enumerate}
If there exists a curve segment $\Gamma$ in $M\setminus M^0$ and a proper rigid motion $\rho\colon \R^3\to\R^3$ such that $f=\rho\circ \tilde{f}$ on $\Gamma$, then $f=\rho\circ \tilde{f}$ on $M$.
\end{thm}

To prove this result, we will again replace $\tilde f$ with $\rho\circ\tilde f$, so that $f=\tilde f$ on $\Gamma$, and show that $f=\tilde f$ on $M$. Let $M^+:=M\setminus M^0$, and $\ol{M^+}$ be the closure of $M^+$ in $M$. Since $M^+$ is connected, Theorem \ref{thm:main} yields that $f=\tilde f$ on $M^+$, and therefore  on $\ol{M^+}$ by continuity. In particular it follows that the second fundamental forms of $f$ and $\tilde f$ agree on $\ol M^+$. Further since $M^0\subset\inte(M)$, we have $\d M^0\subset \ol{M^+}$. Now the following result, which requires conditions (i) and (ii) above, immediately completes the proof of Theorem of \ref{thm:main2}, via the fundamental theorem of surfaces.

\begin{lem}[Sacksteder \cite{sacksteder1962}, Thm. I]\label{lem:sacksteder}
If the second fundamental forms of $f$ and $\tilde f$ agree on $\d M^0$,  then they agree on $M^0$.
\end{lem}

Since Sacksteder's argument is somewhat involved, we  include here a short simple proof of Theorem \ref{thm:main2} for the case where $M^0$ is compact (in which case conditions (i) and (ii) are automatically satisfied). Recall that we just need to check that $f=\tilde f$ on $\inte(M^0)$.  
To this end, let $M^F\subset \inte(M^0)$ be the set of points with a \emph{flat neighborhood}, i.e., a neighborhood mapped by $f$ into a plane. The following fact from Spivak \cite{spivak:v3} is implicit in the works of Hartman-Nirenberg \cite{hartman&nirenberg} and Massey \cite{Mas62} on $\C^2$ surfaces of zero curvature. See also Pogorelov \cite[p. 609]{pogorelov:book} or \cite[p. 79]{sabitov2008} for an extension of this fact to the $\C^1$ category.

\begin{lem}[\cite{spivak:v3}, Cor. 8, p. 243]\label{lem:spivak}
Through every point $p\in\inte(M^0)\setminus M^F$  there passes a curve  $\ell$ with end point(s) on $\d M^0$ such that $f$ maps $\ell$ homeomorphically into a straight line segment or ray in $\R^3$.
\end{lem}

When $M^0$ is compact, $\ell$ has finite length, and so it has two end points $q_i$. Further, since $q_i\in\d M^0\subset \ol{M^+}$, $f(q_i)=\tilde f(q_i)$. So $\tilde f(\ell)$ is a curve joining $f(q_i)$, with the same arc length as $f(\ell)$ by isometry. Consequently $f=\tilde f$ on $\ell$, which yields $f=\tilde f$ on $\inte(M^0)\setminus M^F$.
Next let $C$ be a component of $M^F$. Then $\d C\subset  \ol M^+\cup (M^0\setminus M^F)$. So $f=\tilde f$ on $\d C$, which yields that $f=\tilde f$ on $C$. Indeed,  through each point $p\in C$ there passes  a curve $\ell$ with end points $q_i\in\d C$ such that $f(\ell)$ is a line segment (let $L$ be a complete line passing through $f(p)$ in the plane of $f(C)$,  and $\ell$ be the closure of the component of $f^{-1}(L)\cap C$ containing $p$). Hence, again $f=\tilde f$ on $\ell$, since $f(q_i)=\tilde f(q_i)$. So $f=\tilde f$ on $M^F$, and consequently on $M^0$. 

\begin{note}[The case of zero curvature]
The above argument shows that if $M$ is compact, and $f$, $\tilde f\colon M\to\R^3$ are $\C^2$ isometric immersions with everywhere vanishing curvature, then $f=\tilde f$ on $M$ whenever $f=\tilde f$ on $\d M$. In other words, compact developable surfaces in $\R^3$ are rigid relative to their boundary. 
\end{note}

\appendix

\section{A Short Proof of Cohn-Vossen's Theorem}\label{appendix:cohn-vossen}
Cohn-Vossen proved the first version of his rigidity result in 1927 for positively curved analytic surfaces \cite{cohn-vossen1927},  before extending it to $\C^3$ surfaces in 1936 \cite{cohn-vossen1936}, see Hopf \cite[p. 168]{hopf:book}. The proof of this theorem included in various texts, e.g., \cite{spivak:v5,han:book, chern:book}, is the 1943 argument by Herglotz \cite{herglotz}  based on his celebrated integral formula. Later Wintner \cite{wintner} established the theorem for $\C^2$ surfaces, and Sacksteder \cite{sacksteder1962} extended it to nonnegative curvature; see also \cite[Sec. 6.3]{dajczer} and \cite{guan-shen}. Following the same outline as in the proof of Theorem \ref{thm:main}, we present a proof of Cohn-Vossen's theorem which is even shorter than Herglotz's and works immediately in the $\C^{2,1}$ category.

Let $f$, $\tilde f\colon\S^2\to\R^3$ be $\C^{2,1}$ positively curved isometric immersions, with principal curvatures $k_1\leq k_2$, $\tilde k_1\leq\tilde k_2$  respectively. By the invariance of Gauss curvature, $k_1k_2\equiv\tilde k_1\tilde k_2$. So if, at some point, $k_2<\tilde k_2$, then $\tilde k_1<k_1$, which in turn yields that $\tilde k_1<\tilde k_2$. As is well-known, it is impossible for the last inequality to hold everywhere; because then the principal directions of $\tilde f$ corresponding to $\tilde k_2$, would generate a line field on $\S^2$, in violation of the Poincar\'{e}-Hopf index theorem \cite[Thm. 20, p. 223]{spivak:v3}. So we conclude that $k_2(p)=\tilde k_2(p)$ for some point $p\in\S^2$, which in turn yields that $k_1(p)=\tilde k_1(p)$. Consequently, after a rigid motion, we may assume that $f$ and $\tilde f$ have contact of order $2$ at $p$.

Let $\Omega\subset\R^2$ be an open disk, and $\theta\colon\ol\Omega\to\S^2$ be a smooth map with $\theta(\d\Omega)=p$ such that $\theta\colon\Omega\to\S^2\setminus\{p\}$ is a diffeomorphism. Replace $f$, $\tilde f$ by $f\circ\theta$, $\tilde f\circ\theta$ respectively. Further, as in Section \ref{sec:unique}, set $u:=\l f,e\r$,  $\tilde u:=\l \tilde f,e\r$ for a unit vector $e\in\R^3$. Then $u$, $\tilde u\in \C^{2,1}(\ol\Omega)$,  and $u=\tilde u$ on $\d\Omega$. Again, as in Section \ref{sec:unique}, set $\phi:=u-\tilde u$. For an open set of unit vectors $e$, note that $\phi$ satisfies the linear elliptic equation \eqref{eq:phi} near $\partial \Omega$, which can be put in the  form \eqref{eq:abc} with $c\equiv 0$. Now, since $\phi=0$ on $\d\Omega$,  the unique continuation principle we described in Section \ref{sec:unique} yields that $\phi=0 $ near $\partial\Omega$, which in turn yields that $f=\tilde f$ near $p$. So, by Theorem \ref{thm:main}, $f=\tilde f$ everywhere.

The above argument might also work for nonnegative curvature via an appropriate version of the maximum principle for degenerate equations. In closing, we should recall that Pogorelov \cite{pogorelov:rigidity, pogorelov:book} generalized Cohn-Vossen's theorem to all closed convex surfaces regardless of their regularity in 1952, although that proof remains long and intricate. See also Volkov \cite[Sec. 12.1]{alexandrov:polyhedra} for another approach to Pogorelov's theorem via uniform rigidity estimates.

\section*{Acknowledgments}
We thank Igor Belegradek, Robert Bryant, Jeanne Clelland, Robert Greene, Idjad Sabitov, Andrzej Swiech, and Deane Yang for useful communications. Jellett's work \cite{jellett1849} and its citation in \cite{kagan:book} were first pointed out to us by Sabitov. Further we are grateful to Sabitov for suggesting a simplification of our original proof of Theorem \ref{thm:main2} in the compact case.
\bibliographystyle{abbrv}
\bibliography{references}

\end{document}